\documentclass[reqno,11pt]{amsart}
\usepackage{enumerate,amssymb,amsthm,hyperref}
\usepackage{amsmath,amscd}
\usepackage{mathtools}
\usepackage[nospace,noadjust]{cite}
\usepackage[margin=1in]{geometry}  
\usepackage{tikz}
\usetikzlibrary{decorations.markings,arrows,arrows.meta,positioning}
\definecolor{mygreen}{RGB}{0,134,39}

\tikzset{->-/.style={decoration={
  markings,
  mark=at position #1 with {\arrow{Computer Modern Rightarrow[length=5pt,width=5pt]}}},postaction={decorate}}}
\tikzset{->-rev/.style={decoration={
  markings,
  mark=at position #1 with {\arrow{Computer Modern Rightarrow[length=5pt,width=5pt,reversed]}}},postaction={decorate}}}

\newtheorem{cor}{Corollary}
\newtheorem{thm}[cor]{Theorem}
\newtheorem{prop}[cor]{Proposition}

\theoremstyle{definition}
\newtheorem{rem}[cor]{Remark}
\newtheorem{que}[cor]{Question}
\newtheorem{exa}[cor]{Example}

\newcommand{\prodd}{\operatorname{prod}}
\newcommand{\ov}{\widebar} 
\newcommand{\CA}{\operatorname{CA}}
\newcommand{\Com}{\operatorname{Comm}}
\newcommand{\im}{\operatorname{Im}}
\newcommand{\Ker}{\operatorname{Ker}}
\newcommand{\PMod}{\mathcal{PM}}
\newcommand{\M}{\mathcal M}
\renewcommand{\P}{\mathcal P}

\renewcommand{\S}{\mathfrak S}
\newcommand{\Z}{\mathbb Z}
\newcommand{\Aut}{\operatorname{Aut}}
\newcommand{\Out}{\operatorname{Out}}
\newcommand{\Inn}{\operatorname{Inn}}

  \makeatletter
  \newcommand*\if@single[3]{%
    \setbox0\hbox{${\mathaccent"0362{#1}}^H$}%
    \setbox2\hbox{${\mathaccent"0362{\kern0pt#1}}^H$}%
    \ifdim\ht0=\ht2 #3\else #2\fi
    }
  \newcommand*\rel@kern[1]{\kern#1\dimexpr\macc@kerna}
  \newcommand*\widebar[1]{\@ifnextchar^{{\wide@bar{#1}{0}}}{\wide@bar{#1}{1}}}
  \newcommand*\wide@bar[2]{\if@single{#1}{\wide@bar@{#1}{#2}{1}}{\wide@bar@{#1}{#2}{2}}}
  \newcommand*\wide@bar@[3]{%
    \begingroup
    \def\mathaccent##1##2{%
      \if#32 \let\macc@nucleus\first@char \fi
      \setbox\z@\hbox{$\macc@style{\macc@nucleus}_{}$}%
      \setbox\tw@\hbox{$\macc@style{\macc@nucleus}{}_{}$}%
      \dimen@\wd\tw@
      \advance\dimen@-\wd\z@
      \divide\dimen@ 3
      \@tempdima\wd\tw@
      \advance\@tempdima-\scriptspace
      \divide\@tempdima 10
      \advance\dimen@-\@tempdima
      \ifdim\dimen@>\z@ \dimen@0pt\fi
      \rel@kern{0.6}\kern-\dimen@
      \if#31
        \overline{\rel@kern{-0.9}\kern\dimen@\macc@nucleus\rel@kern{0.4}\kern\dimen@}%
        \advance\dimen@0.4\dimexpr\macc@kerna
        \let\final@kern#2%
        \ifdim\dimen@<\z@ \let\final@kern1\fi
        \if\final@kern1 \kern-\dimen@\fi
      \else
        \overline{\rel@kern{-0.9}\kern\dimen@#1}%
      \fi
    }%
    \macc@depth\@ne
    \let\math@bgroup\@empty \let\math@egroup\macc@set@skewchar
    \mathsurround\z@ \frozen@everymath{\mathgroup\macc@group\relax}%
    \macc@set@skewchar\relax
    \let\mathaccentV\macc@nested@a
    \if#31
      \macc@nested@a\relax111{#1}%
    \else
      \def\gobble@till@marker##1\endmarker{}%
      \futurelet\first@char\gobble@till@marker#1\endmarker
      \ifcat\noexpand\first@char A\else
        \def\first@char{}%
      \fi
      \macc@nested@a\relax111{\first@char}%
    \fi
    \endgroup
  }
  \makeatother

\makeatletter
\@namedef{subjclassname@2020}{\textup{2020} Mathematics Subject Classification}
\makeatother

\begin{document}

\date{\today}
\title{Artin groups of types $F_4$ and $H_4$ are not commensurable with that of type $D_4$}

\author{Ignat Soroko}
\address{Department of Mathematics\\
Florida State University\\
Tallahassee\\ FL 32306\\ USA}
\email{ignat.soroko@gmail.com}



\begin{abstract}
In a recent article, Cumplido and Paris studied the question of commensurability between Artin groups of spherical type. Their analysis left six cases undecided, for the following pairs of Artin groups: $(F_4,D_4)$, $(H_4,D_4)$, $(F_4,H_4)$, $(E_6,D_6)$, $(E_7,D_7)$, and $(E_8,D_8)$. In this note we resolve the first two of these cases, namely, we show that the Artin groups of types $F_4$ and $H_4$ are not commensurable with that of type $D_4$. As a key step, we realize the abstract commensurator of the Artin group of type $D_4$ as the extended mapping class group of the torus with three punctures. We also find the automorphism group of the Artin group of type $D_4$ and obtain a description of torsion elements, their orders and conjugacy classes in all irreducible Artin groups of spherical type modulo their centers.
\end{abstract}

\maketitle

\section{Introduction}
Artin groups (also called the Artin--Tits groups) are a natural generalization of the classical braid groups. They admit a presentation in which the only defining relations are (generalized) braid relations between some pairs of generators, which is usually encoded via a graph whose vertices are generators and edges are labeled by $\{3,\dots,\infty\}$ (the associated Coxeter graph, see Section~\ref{sec:2} for definitions). Among all Artin groups, one subclass stands out for richness of their structure and intrinsic beauty, the class of Artin groups of spherical type. They are exactly the Artin groups whose associated Coxeter group is finite. The latter were described in terms of their defining graphs by Coxeter himself, who completely classified all connected graphs corresponding to finite Coxeter groups. The answer is a list of four infinite series $A_n$, $B_n$, $D_n$, $I_2(m)$ and six standalone groups $E_6$, $E_7$, $E_8$, $F_4$, $H_3$, and $H_4$, presented in Figure~\ref{fig:sph}. This list also describes all directly indecomposable Artin groups of spherical type up to isomorphism, as was shown in~\cite{Paris1}. One can also ask which Artin groups of spherical type are commensurable, i.e.\ allow embeddings into one another `up to finite index'. This question was studied in the article of Cumplido and Paris~\cite{CP}, 
where the authors show that to be commensurable two Artin groups of spherical type must have the same ranks, and their irreducible components should be pairwise commensurable. Cumplido and Paris completely resolve the question of commensurability when one of the groups is of type $A_n$, and their result states that for each rank $n\ge3$, the only spherical Artin group commensurable with the Artin group of type $A_n$ is that of type $B_n$, whereas for rank $n=2$ all Artin groups of type $I_2(m)$, $m\ge3$, are commensurable with $A_2$. To obtain the complete classification of Artin groups of spherical type up to commensurability, six remaining cases need to be resolved, for Artin groups corresponding to the pairs: $(F_4,D_4)$, $(H_4,D_4)$, $(F_4,H_4)$, $(E_6,D_6)$, $(E_7,D_7)$, and $(E_8,D_8)$. In this note we resolve the first two out of these six cases. Namely, we prove:\let\thefootnote\relax\footnotetext{\copyright{} 2021. This manuscript version is made available under the CC--BY--NC--ND 4.0 license\\ \url{http://creativecommons.org/licenses/by-nc-nd/4.0/}}
\begin{thm}
Artin groups of types $F_4$ and $H_4$ are not commensurable with the Artin group of type~$D_4$. 
\end{thm}
The proof of this theorem is given in Propositions~\ref{prop:d4h4} and~\ref{prop:d4f4} below. 
Note that we still do not know if Artin groups of types $F_4$ and $H_4$ are commensurable with each other.

Let $A[\Gamma]$ denote the Artin group of a spherical type $\Gamma$, $\ov{A[\Gamma]}$ denote its quotient by the center, and $\Com(G)$ be the abstract commensurator of a group $G$ (see Section~\ref{sec:2} for definitions).

To establish the non-commensurability of $A[\Gamma]$ with $A[A_n]$ Cumplido and Paris implement the following strategy. An assumption that $A[\Gamma]$ and $A[A_n]$ are commensurable implies that there exists an embedding $\Phi$ of $\ov{A[\Gamma]}$ into the abstract commensurator $\Com(\ov{A[A_n]})$ of $\ov{A[A_n]}$, which is known to be isomorphic to the extended mapping class group of the sphere with punctures. To prove that such an embedding $\Phi$ does not exist, Cumplido and Paris analyze its projections onto the finite quotient of $\Com(\ov{A[A_n]})$ by a certain bi-orderable group (isomorphic to the pure braid group). Assisted by computer, they classify all such projections, and in most of the cases they are able to exhibit a generalized torsion element in the (bi-orderable) kernel of the projection, which gives a contradiction. The remaining `hard' cases are dealt with by finding an element in the image of $\Phi$ which is pseudo-Anosov but whose centralizer is not virtually cyclic, thus yielding a contradiction as well.

Our approach is similar: we use a realization of $\ov{A[D_4]}$ as the pure mapping class group $\PMod(\Sigma_{1,0},\P_3)$ of the torus with three punctures~\cite{S} (a fact which deserves to be more broadly known in the mathematical community). This yields a description of the commensurator $\Com(\ov{A[D_4]})$ of $\ov{A[D_4]}$ as the extended mapping class group $\M^*(\Sigma_{1,0},\P_3)$ in Theorem~\ref{thm:comd4}. As in~\cite{CP}, we need to prove that embeddings of $\ov{A[H_4]}$ and $\ov{A[F_4]}$ into $\Com(\ov{A[D_4]})$ do not exist. To show this for $\ov{A[H_4]}$, we look at the torsion elements. It turns out that there exist elements of order $15$ in $\ov{A[H_4]}$ but not in $\Com(\ov{A[D_4]})$, which yields a contradiction. To prove that an embedding of $\ov{A[F_4]}$ to $\Com(\ov{A[D_4]})$ does not exist, we classify all epimorphisms from $\ov{A[F_4]}$ onto the finite quotient of $\Com(\ov{A[D_4]})$ by the (bi-orderable) pure Artin group of type $D_4$ modulo center.  
The target finite group is quite large (it has order~$1156$), so we use the computer algebra system \textsc{Magma}~\cite{Magma} for that, which has a functionality (due to Derek Holt) for computing all homomorphisms, including non-surjective ones, from a finitely presented group to a finite group, up to the conjugacy in the latter. In all cases but one we are able to exhibit generalized torsion elements in the bi-orderable kernel, which gives a contradiction. The remaining `hard' case appears in the analogous situation as in~\cite{CP}, but in the case of the punctured torus we do not have a handy theorem that allows us to use properties of pseudo-Anosov elements. Instead, to obtain a contradiction, we utilize the symmetry of the Coxeter graph of $F_4$, together with the result of Behrstock and Margalit~\cite{BM} which states that every injection of a finite index subgroup of the extended mapping class group $\M^*(\Sigma_{1,0},\P_3)$ into $\M^*(\Sigma_{1,0},\P_3)$ is induced by a conjugation in $\M^*(\Sigma_{1,0},\P_3)$.

Our article is organized as follows. Section~\ref{sec:2} contains the necessary definitions and key auxiliary results that we use. In Section~\ref{sec:3}, by realizing $\ov{A[D_4]}$ as the extended mapping class group of the torus with three punctures $(\Sigma_{1,0},\P_3)$, we describe the abstract commensurator $\Com(\ov{A[D_4]})$, and answer a question of Crisp and Paris from~\cite{CrispParisBD} about the structure of the automorphism group and the outer automorphism group for $A[D_4]$ and $\ov{A[D_4]}$. Namely, we establish in Theorem~\ref{thm:comd4} and in Corollary~\ref{cor:aut} that:
\begin{gather*}
\Com(\ov{A[D_4]})\cong\Aut(\ov{A[D_4]})\cong\Aut(A[D_4])\cong\ov{A[D_4]}\rtimes (\S_3\times \Z_2)\cong \M^*(\Sigma_{1,0},\P_3),\\
\Out(\ov{A[D_4]})\cong\Out(A[D_4])\cong\S_3\times \Z_2,
\end{gather*}
where $\S_3$ is the symmetric group on three elements and $\Z_2$ is the cyclic group of order two. In Section~\ref{sec:4} we describe the orders of torsion elements in all Artin groups of spherical type modulo their centers, and give representatives of the torsion elements of each order up to the conjugacy (Theorem~\ref{th:orders} and Corollary~\ref{cor:phi}). Using this information, we establish in Section~\ref{sec:5} that the Artin groups $A[H_4]$ and $A[D_4]$ are not commensurable (Proposition~\ref{prop:d4h4}). 
Using the strategy outlined above, we establish in Section~\ref{sec:6} that the Artin groups $A[F_4]$ and $A[D_4]$ are not commensurable (Proposition~\ref{prop:d4f4}). 

In~\cite{Waj} Wajnryb asks which Artin groups admit embeddings into mapping class groups of compact surfaces (with or without boundary). Since all boundary Dehn twists are central in the mapping class group, it is natural to ask the following related question for the central quotients~$\ov{A[\Gamma]}$:
\begin{que}\label{que:emb}
For which Artin groups $A[\Gamma]$ their central quotients $\ov{A[\Gamma]}$ admit embeddings into mapping class groups of punctured surfaces?
\end{que}
Our results give partial answer to this question. Namely, Corollary~\ref{cor:h4mcg} establishes that the group ${}\ov{A[H_4]}$ does not admit any embedding into the extended mapping class group $\M^*(\Sigma_{1,0},\P_3)$ of the three times punctured torus, and Proposition~\ref{prop:d4f4} implies that there does not exist an embedding of ${}\ov{A[F_4]}$ into $\M^*(\Sigma_{1,0},\P_3)$ such that the image has finite index in the latter. In the Example~\ref{ex:13} we analyze a natural `geometric' candidate for an embedding of $\ov{A[F_4]}$ into $\M^*(\Sigma_{1,0},\P_3)$, which sends the standard generators of $\ov{A[F_4]}$ to Dehn twists and half-twists, and show that (unfortunately) it is not injective.

\subsection*{Acknowledgments} The author is grateful to Mar\'ia Cumplido and Luis Paris for useful conversations about their work~\cite{CP}, and to Aaron Calderon for his catalyzing remarks which led the author to an idea of writing this paper. The author thanks the organizers of the `young Geometric Group Theory VIII' conference in Bilbao, June 30--July 5, 2019, and the `Virtual Geometric Group Theory' conference at CIRM, Luminy, June 1--5, 2020, which allowed for these conversations to happen. The author extends his gratitude to Sang-Jin Lee for answering questions on the paper~\cite{LL} and to Derek Holt and Alexander Hulpke for their helpful consultations on the computer algebra systems \textsc{Magma} and \textsf{GAP}, respectively. 
The author also thanks the referee for careful reading of the manuscript and providing useful suggestions which improved the quality of this paper. This work was partially supported by the AMS--Simons Travel Grant.

\section{Preliminaries}\label{sec:2}
In this section we collect the relevant definitions and key auxiliary results used in this work.

\begin{figure}
\begin{center}
\begin{tikzpicture}[scale=0.7]
\begin{scope}[xshift=-5cm,yshift=-3.5cm,scale=1.33] 
\fill (0,0) circle (2.3pt) node [below=2pt] {\small$s_1$};
\fill (1,0) circle (2.3pt) node [below=2pt] {\small$s_2$};
\fill (2,0) circle (2.3pt) node [below=2pt] {\small$s_3$};
\fill (4,0) circle (2.3pt) node [below=2pt] {\small$s_{n-1}$};
\fill (5,0) circle (2.3pt) node [below=2pt] {\small$s_n$};
\draw [thick] (0,0)--(1,0)--(2,0)--(2.5,0);
\draw [thick,dashed] (2.5,0)--(3.5,0);
\draw [thick] (3.5,0)--(4,0)--(5,0);
\draw (-1.8,0) node {$A_n$, $(n\ge1)$:};	
\end{scope}

\begin{scope}[xshift=-5cm,yshift=-5.5cm,scale=1.33] 
\fill (0,0) circle (2.3pt) node [below=2pt] {\small$s_1$};
\fill (1,0) circle (2.3pt) node [below=2pt] {\small$s_2$};
\fill (2,0) circle (2.3pt) node [below=2pt] {\small$s_3$};
\fill (4,0) circle (2.3pt) node [below=2pt] {\small$s_{n-1}$};
\fill (5,0) circle (2.3pt) node [below=2pt] {\small$s_n$};
\draw [thick] (0,0)--(1,0)--(2,0)--(2.5,0);
\draw [thick,dashed] (2.5,0)--(3.5,0);
\draw [thick] (3.5,0)--(4,0)--(5,0);
\draw (0.5,0.3) node {\small4};
\draw (-1.8,0) node {$B_n$, $(n\ge2)$:};	
\end{scope}

\begin{scope}[xshift=-5cm,yshift=-7.75cm,scale=1.33] 
\fill (0,0.5) circle (2.3pt) node [below=2pt] {\small$s_1$};
\fill (0,-0.5) circle (2.3pt) node [below=2pt] {\small$s_2$};
\fill (1,0) circle (2.3pt) node [below=2pt] {\small$s_3$};
\fill (2,0) circle (2.3pt) node [below=2pt] {\small$s_4$};
\fill (4,0) circle (2.3pt) node [below=2pt] {\small$s_{n-1}$};
\fill (5,0) circle (2.3pt) node [below=2pt] {\small$s_n$};
\draw [thick] (0,0.5)--(1,0)--(0,-0.5);
\draw [thick] (1,0)--(2,0)--(2.5,0);
\draw [thick,dashed] (2.5,0)--(3.5,0);
\draw [thick] (3.5,0)--(4,0)--(5,0);
\draw (-1.8,0) node {$D_n$, $(n\ge4)$:};	
\end{scope}

\begin{scope}[xshift=5.5cm,yshift=-3.5cm,scale=1.33] 
\fill (0,0) circle (2.3pt) node [below=2pt] {\small$s_1$};
\fill (1,0) circle (2.3pt) node [below=2pt] {\small$s_3$};
\fill (2,0) circle (2.3pt);\draw (1.75,-0.35) node {\small$s_4$};
\fill (3,0) circle (2.3pt) node [below=2pt] {\small$s_5$};
\fill (4,0) circle (2.3pt) node [below=2pt] {\small$s_6$};
\fill (2,-1) circle (2.3pt) node [left=2pt] {\small$s_2$};
\draw [thick] (0,0)--(1,0)--(2,0)--(3,0)--(4,0);
\draw [thick] (2,0)--(2,-1);
\draw (-1,0) node {$E_6$:};	
\end{scope}

\begin{scope}[xshift=5.5cm,yshift=-6cm,scale=1.33] 
\fill (0,0) circle (2.3pt) node [below=2pt] {\small$s_1$};
\fill (1,0) circle (2.3pt) node [below=2pt] {\small$s_3$};
\fill (2,0) circle (2.3pt);\draw (1.75,-0.35) node {\small$s_4$};
\fill (3,0) circle (2.3pt) node [below=2pt] {\small$s_5$};
\fill (4,0) circle (2.3pt) node [below=2pt] {\small$s_6$};
\fill (5,0) circle (2.3pt) node [below=2pt] {\small$s_7$};
\fill (2,-1) circle (2.3pt) node [left=2pt] {\small$s_2$};
\draw [thick] (0,0)--(1,0)--(2,0)--(3,0)--(4,0)--(5,0);
\draw [thick] (2,0)--(2,-1);
\draw (-1,0) node {$E_7$:};	
\end{scope}

\begin{scope}[xshift=5.5cm,yshift=-8.5cm,scale=1.33] 
\fill (0,0) circle (2.3pt) node [below=2pt] {\small$s_1$};
\fill (1,0) circle (2.3pt) node [below=2pt] {\small$s_3$};
\fill (2,0) circle (2.3pt);\draw (1.75,-0.35) node {\small$s_4$};
\fill (3,0) circle (2.3pt) node [below=2pt] {\small$s_5$};
\fill (4,0) circle (2.3pt) node [below=2pt] {\small$s_6$};
\fill (5,0) circle (2.3pt) node [below=2pt] {\small$s_7$};
\fill (6,0) circle (2.3pt) node [below=2pt] {\small$s_8$};
\fill (2,-1) circle (2.3pt) node [left=2pt] {\small$s_2$};
\draw [thick] (0,0)--(1,0)--(2,0)--(3,0)--(4,0)--(5,0)--(6,0);
\draw [thick] (2,0)--(2,-1);
\draw (-1,0) node {$E_8$:};	
\end{scope}

\begin{scope}[xshift=-5cm,yshift=-10.5cm,scale=1.33] 
\fill (0,0) circle (2.3pt) node [below=2pt] {\small$s_1$};
\fill (1,0) circle (2.3pt) node [below=2pt] {\small$s_2$};
\fill (2,0) circle (2.3pt) node [below=2pt] {\small$s_3$};
\fill (3,0) circle (2.3pt) node [below=2pt] {\small$s_4$};
\draw [thick] (0,0)--(1,0)--(2,0)--(3,0);
\draw (1.5,0.3) node {\small4};
\draw (-1,0) node {$F_4$:};	
\end{scope}

\begin{scope}[xshift=-5cm,yshift=-12.5cm,scale=1.33] 
\fill (0,0) circle (2.3pt) node [below=2pt] {\small$s_1$};
\fill (1,0) circle (2.3pt) node [below=2pt] {\small$s_2$};
\fill (2,0) circle (2.3pt) node [below=2pt] {\small$s_3$};
\fill (3,0) circle (2.3pt) node [below=2pt] {\small$s_4$};
\draw [thick] (0,0)--(1,0)--(2,0)--(3,0);
\draw (0.5,0.3) node {\small5};
\draw (-1,0) node {$H_4$:};	
\end{scope}

\begin{scope}[xshift=2.75cm,yshift=-12.5cm,scale=1.33] 
\fill (0,0) circle (2.3pt) node [below=2pt] {\small$s_1$};
\fill (1,0) circle (2.3pt) node [below=2pt] {\small$s_2$};
\fill (2,0) circle (2.3pt) node [below=2pt] {\small$s_3$};
\draw [thick] (0,0)--(1,0)--(2,0);
\draw (0.5,0.3) node {\small5};
\draw (-1,0) node {$H_3$:};	
\end{scope}

\begin{scope}[xshift=11.75cm,yshift=-12.5cm,scale=1.33] 
\fill (0,0) circle (2.3pt) node [below=2pt] {\small$s_1$};
\fill (1,0) circle (2.3pt) node [below=2pt] {\small$s_2$};
\draw [thick] (0,0)--(1,0);
\draw (0.5,0.3) node {\small $m$};
\draw(-2.2,0.25) node {$I_2(m),$};
\draw(-2.2,-0.25) node {\small$(m\ge5,\, m\ne\infty)$};
\draw(-0.65,0) node {:};
\end{scope}

\end{tikzpicture}
\end{center}
\caption{Coxeter graphs of irreducible Artin groups of spherical type.\label{fig:sph}}
\end{figure}
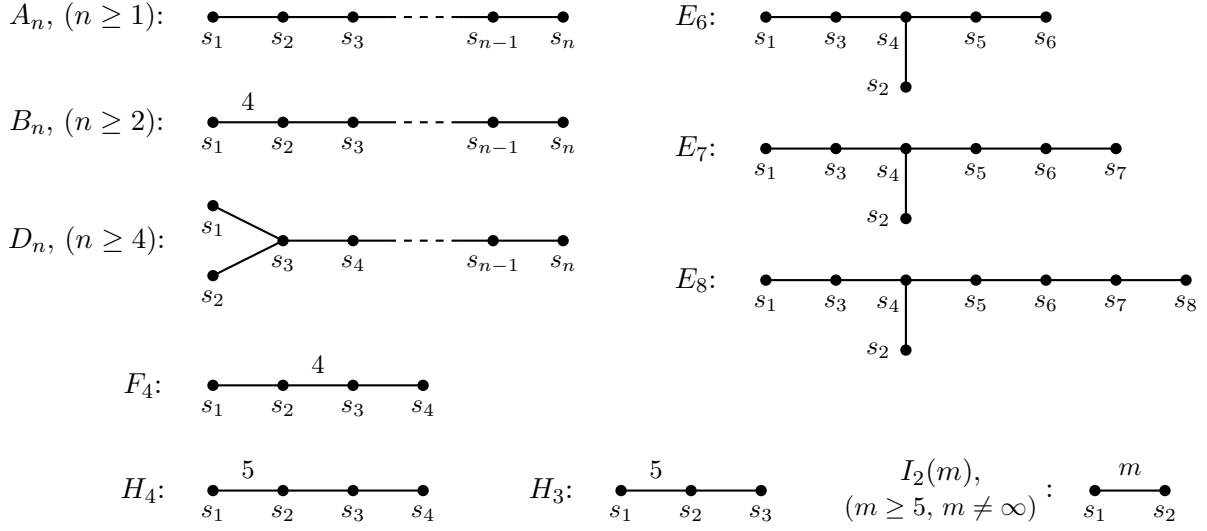

Recall that a \textit{Coxeter matrix} over a finite set $S$ is a symmetric matrix $(m_{st})_{s,t\in S}$ with entries in $\{1,2,\dots,\infty\}$, such that $m_{ss}=1$ for all $s\in S$ and $m_{st}=m_{ts}\ge2$ if $s\ne t$. A Coxeter matrix can be encoded by the corresponding \textit{Coxeter graph} $\Gamma$ having $S$ as the set of vertices. Two distinct vertices $s,t\in S$ are connected with an edge in $\Gamma$ if $m_{st}\ge 3$, and this edge is labeled with $m_{st}$ if $m_{st}\ge 4$. The \textit{Artin group} associated to $\Gamma$ is the group $A[\Gamma]$ given by the presentation:
\[
A[\Gamma]=\langle S\mid \prodd(s,t,m_{st})=\prodd(t,s,m_{ts}), \text{ for all } s\ne t,\,m_{st}\ne\infty\rangle,
\]
where $\prodd(s,t,m_{st})$ is the word $stst\dots$ of length $m_{st}\ge2$. The \textit{Coxeter group} $W[\Gamma]$ of $\Gamma$ is the quotient of $A[\Gamma]$ by all relations of the form $s^2=1$, $s\in S$. The Artin group $A[\Gamma]$ is \textit{of spherical type} if the corresponding Coxeter group $W[\Gamma]$ is finite. The group $A[\Gamma]$ is called \textit{irreducible} if the corresponding Coxeter graph $\Gamma$ is connected. All Artin groups of spherical type are classified up to isomorphism in~\cite{Paris1}: two such groups are isomorphic if and only if their Coxeter graphs are isomorphic (as edge-labeled graphs). The list of all Coxeter graphs corresponding to irreducible Artin groups of spherical type is given in Figure~\ref{fig:sph}. (Note that for small values of $n,m$ we have the following isomorphisms of labeled graphs: $B_1\cong A_1$, $D_3\cong A_3$, $I_2(3)\cong A_2$, $I_2(4)\cong B_2$, so we specified inequalities on $n,m$ to obtain a duplicate-free list.)

It is known~\cite{BrieskornSaito} that the center of any irreducible Artin group of spherical type is an infinite cyclic group. We denote $\ov{A[\Gamma]}$ the quotient of the Artin group $A[\Gamma]$ by its center. We also denote by $\CA[\Gamma]$ the kernel of the natural projection $A[\Gamma]\longrightarrow W[\Gamma]$, which is called the \textit{pure Artin group} corresponding to $\Gamma$. Denote $\ov{\CA[\Gamma]}$ the quotient of $\CA[\Gamma]$ by its center $Z(\CA[\Gamma])$. In the following Proposition we collect some properties of the groups $\CA[\Gamma]$ and $\ov{\CA[\Gamma]}$ that we will use:
\begin{prop}\label{prop:ca}
Let $A[\Gamma]$ be an irreducible Artin group of spherical type. Then:
\pushQED{\qed}
\begin{enumerate}
\item The center $Z(\CA[\Gamma])$ of $\CA[\Gamma]$ is an infinite cyclic group~\cite[Cor.\,7]{CP};
\item $\ov{\CA[\Gamma]}$ has trivial center and embeds into $\ov{A[\Gamma]}$ as a subgroup of finite index~\cite[Cor.\,7]{CP};
\item $\CA[\Gamma]\cong \ov{\CA[\Gamma]}\times Z(\CA[\Gamma])$~\cite[Prop.\,6\,(2)]{CP};
\item $\CA[\Gamma]$ (and hence $\ov{\CA[\Gamma]}$) is bi-orderable~\cite[Th.\,3]{Mar}.\qedhere
\end{enumerate}
\popQED
\end{prop}

Let $G$ be a group. Recall that an element $\alpha\in G$ is called a \textit{generalized torsion element} if there exists $n\ge1$ and elements $\beta_1,\dots,\beta_n\in G$ such that $(\beta_1\alpha\beta_1^{-1})(\beta_2\alpha\beta_2^{-1})\dots(\beta_n\alpha\beta_n^{-1})=1$. It is known that if a group $G$ is bi-orderable it cannot have generalized torsion elements~\cite{RZ}. In particular, groups $\CA[\Gamma]$ for irreducible Artin groups of spherical type do not contain generalized torsion elements, by Proposition~\ref{prop:ca}.

Recall that two groups $G_1$ and $G_2$ are called \textit{commensurable} if there exist finite index subgroups $H_1\subseteq G_1$, $H_2\subseteq G_2$ such that $H_1$ is isomorphic to $H_2$. For each group $G$ there is a canonically constructed group $\Com(G)$ called the (abstract) \textit{commensurator} of $G$, which is defined as follows. Let $\widetilde\Com(G)$ be the set of all triples $(U,V,f)$ where $U,V$ are finite index subgroups of $G$ and ${f\colon U\to V}$ an isomorphism. Let $\sim$ be the equivalence relation on $\widetilde\Com(G)$ such that $(U,V,f)\sim(U',V',f')$ if there exists a finite index subgroup $W$ of $U\cap U'$ such that $f|_W=f'|_W$. As a set, the commensurator $\Com(G)$ is equal to $\Com(G)=\widetilde\Com(G)/{\sim}$. If $[f_1],[f_2]\in\Com(G)$ are two equivalence classes corresponding to triples $(U_1,V_1,f_1)$, $(U_2,V_2,f_2)\in\widetilde\Com(G)$, then the composition of classes $[f_1]\circ[f_2]$ is defined as the equivalence class corresponding to the triple $\bigl(f_2^{-1}(U_1\cap V_2),f_1(U_1\cap V_2),(f_1\circ f_2)|_{f_2^{-1}(U_1\cap V_2)}\bigr)$. It can be shown that under this operation $\Com(G)$ forms a group, and if two groups $G_1$ and $G_2$ are commensurable, then their commensurators are isomorphic. A very readable account of basic properties of commensurators is given in the survey of Studenmund~\cite{St}.

Cumplido and Paris established the following important result, which is key to the proofs in~\cite{CP} and in the present article:
\begin{thm}[\protect{\cite[Prop.\,6\,(3),(4)]{CP}}]\label{thm:comm}
Let $A[\Gamma]$ and $A[\Omega]$ be two irreducible Artin groups of spherical type. Then
\pushQED{\qed}
\begin{enumerate}
\item $A[\Gamma]$ and $A[\Omega]$ are commensurable if and only if $\ov{A[\Gamma]}$ and $\ov{A[\Omega]}$ are commensurable;
\item The group $\ov{A[\Gamma]}$ injects into its commensurator $\Com\big(\ov{A[\Gamma]}\big)$.\qedhere
\end{enumerate}
\popQED
\end{thm}

One more ingredient that we will need is the mapping class group of a surface, the definition of which we briefly recall. Let $\Sigma=\Sigma_{g,b}$ be the orientable surface of genus $g$ with $b$ boundary components, and let $\P_n$ be a collection of $n$ different points in the interior of $\Sigma_{g,b}$. The mapping class group $\M(\Sigma_{g,b},\P_n)$ of the pair $(\Sigma_{g,b},\P_n)$ is the group of orientation-preserving homeomorphisms of $\Sigma_{g,b}$, identical on the boundary and permuting the set $\P_n$, considered up to isotopies identical on the boundary and fixing $\P_n$ pointwise. If we allow orientation-reversing homeomorphisms in the above definition, we get the notion of the \textit{extended mapping class group} of the pair $(\Sigma_{g,b},\P_n)$, which is denoted $\M^*(\Sigma_{g,b},\P_n)$. If a surface $\Sigma$ has nonempty boundary, then $\M^*(\Sigma,\P_n)=\M(\Sigma,\P_n)$, otherwise $\M(\Sigma,\P_n)$ is a subgroup of index $2$ in $\M^*(\Sigma,\P_n)$. The \textit{pure mapping class group} of the pair $(\Sigma_{g,b},\P_n)$ is the finite index subgroup $\PMod(\Sigma_{g,b},\P_n)$ of $\M(\Sigma_{g,b},\P_n)$ which fixes the set $\P_n$ pointwise.

\section{The mapping class group of the torus with three punctures}\label{sec:3} 
It is known that the pure mapping class group of the three times punctured torus is isomorphic to $\ov{A[D_4]}$, the Artin group of type $D_4$ modulo its center. This fact is implicit in the description of the general presentations for the mapping class groups of punctured surfaces given by Labru\`ere and Paris in~\cite{LP}, and can also be deduced from the Gervais presentation~\cite{Ger}, as it was shown in~\cite{S}.
We will use it to describe the abstract commensurator of $\ov{A[D_4]}$ in the following theorem.

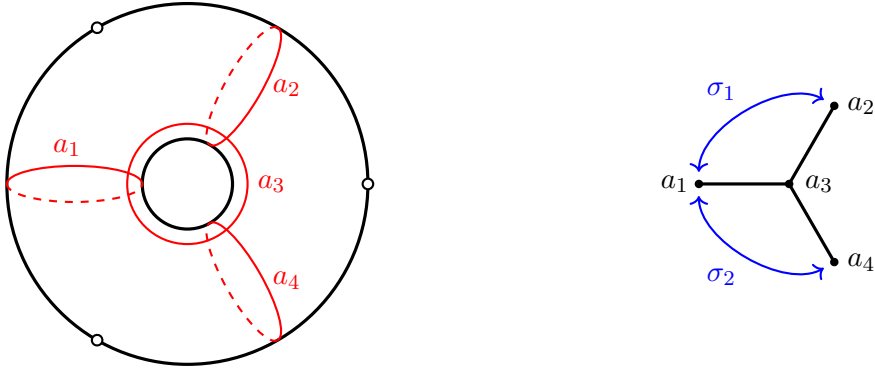
\begin{figure}[ht!]
\begin{center}
\begin{tikzpicture}[scale=0.8]
\draw [very thick] (0,0) circle (0.75cm);
\draw [black,very thick] (0,0) circle (3cm);

\color{red}
\draw[thick] (-3,0) arc (180:0:1.125cm and 0.30cm);
\draw[dashed,thick] (-3,0) arc (180:360:1.125cm and 0.30cm);
\draw[black,thick,fill=white,radius=2.5pt] (3,0) circle;
\begin{scope}[rotate=-120]
\draw[thick] (-3,0) arc (180:0:1.125cm and 0.30cm);
\draw[thick,dashed] (-3,0) arc (180:360:1.125cm and 0.30cm);
\draw[black,thick,fill=white,radius=2.5pt] (3,0) circle;
\end{scope}
\begin{scope}[rotate=120]
\draw[thick,dashed] (-3,0) arc (180:0:1.125cm and 0.30cm);
\draw[thick] (-3,0) arc (180:360:1.125cm and 0.30cm);
\draw [black,thick,fill=white,radius=2.5pt] (3,0) circle;
\end{scope}
\draw (1.65,-1.6) node {$a_4$};
\draw (1.65,1.6) node {$a_2$};
\draw (-2,0.6) node {$a_1$};
\draw [thick] (0,0) circle (1cm);
\draw (1.4,0) node {$a_3$};

\begin{scope}[xshift=10cm]
\color{black}
\fill (0,0) circle (2pt);
\fill (-1.5,0) circle (2pt);
\draw[very thick] (-1.5,0)--(0,0);
\begin{scope}[rotate=120]
\fill (-1.5,0) circle (2pt);
\draw[very thick] (-1.5,0)--(0,0);
\end{scope}
\begin{scope}[rotate=-120]
\fill (-1.5,0) circle (2pt);
\draw[very thick] (-1.5,0)--(0,0);
\end{scope}
\draw (0.5,0) node {$a_3$};
\draw (1.2,-1.3) node {$a_4$};
\draw (1.2,1.3) node {$a_2$};
\draw (-1.9,0) node {$a_1$};
\draw[blue,thick,<->,looseness=0.8] (-1.5,0.2) to [out=90, in=150, edge label=$\sigma_1$] (0.577,1.4);
\begin{scope}[rotate=120]
\draw[blue,thick,<->,looseness=0.8] (-1.5,0.2) to [out=90, in=150, edge label=$\sigma_2$] (0.577,1.4);
\end{scope}
\end{scope}

\end{tikzpicture}
\caption{The Dehn twist generators of $\protect\ov{A[D_4]}$ and the action of\/ $\S_3$ on them.\label{fig:torus}}
\end{center}
\end{figure}

\begin{thm}\label{thm:comd4}
The abstract commensurator $\Com(\ov{A[D_4]})$ of $\ov{A[D_4]}$ is isomorphic to the extended mapping class group $\M^{*}(\Sigma_{1,0},\P_3)$ of the torus with three punctures, which has the following structure:
\[
\Com(\ov{A[D_4]})\cong\M^{*}(\Sigma_{1,0},\P_3)\cong \ov{A[D_4]}\rtimes (\S_3\times \Z_2),
\]
where $\ov{A[D_4]}\cong\PMod(\Sigma_{1,0},\P_3)$ is the pure mapping class group of $(\Sigma_{1,0},\P_3)$, $\S_3$ is the symmetric group of order $6$ consisting of rotations of the torus permuting the punctures, and $\Z_2$ is generated by the orientation-reversing reflection fixing all punctures pointwise (i.e.\ the reflection in the `horizontal' plane of the torus in Figure~\ref{fig:torus}).

In the Dehn twist generators $a_1,a_2,a_3,a_4$ about the curves in Figure~\ref{fig:torus}, the group $\PMod(\Sigma_{1,0},\P_3)$ has a presentation:
\begin{multline*}
\PMod(\Sigma_{1,0},\P_3)=\langle a_1,a_2,a_3,a_4\mid a_1a_3a_1=a_3a_1a_3,\quad a_3a_2a_3=a_2a_3a_2,\quad a_3a_4a_3=a_4a_3a_4, \\
a_1a_2=a_2a_1,\quad a_1a_4=a_4a_1,\quad a_2a_4=a_4a_2,\quad (a_1a_2a_3a_4)^3=1\rangle,
\end{multline*}
which coincides with the standard presentation of\/ $\ov{A[D_4]}$. If we denote a pair of transpositions in\/ $\S_3$ as $\sigma_1,\sigma_2$, and the generator of\/ $\Z_2$ as $\iota$, we will have: $\S_3=\langle \sigma_1,\sigma_2\mid \sigma_1\sigma_2\sigma_1=\sigma_2\sigma_1\sigma_2,\sigma_1^2=\sigma_2^2=1\rangle$, $\Z_2=\langle\iota\mid\iota^2=1\rangle$, and the semidirect product structure on $\M^{*}(\Sigma_{1,0},\P_3)$ can be given by the following conjugation action:
\begin{gather*}
a_1^{\sigma_1}=a_2,\quad a_2^{\sigma_1}=a_1,\quad a_3^{\sigma_1}=a_3,\quad a_4^{\sigma_1}=a_4,\\
a_1^{\sigma_2}=a_4,\quad a_4^{\sigma_2}=a_1,\quad a_2^{\sigma_2}=a_2,\quad a_3^{\sigma_2}=a_3,\\
\sigma_1^\iota=\sigma_1,\quad \sigma_2^\iota=\sigma_2, \quad a_k^\iota = a_k^{-1},\quad \text{for $k=1,2,3,4$}. 
\end{gather*}
\end{thm}

\begin{proof}
It is known that for an arbitrary surface $\Sigma$ with a set of $n$ punctures $\P_n$ (and possibly with boundary), the groups $\PMod(\Sigma,\P_n)$ and $\M(\Sigma,\P_n)$ fit into the short exact sequence:
\[
1\longrightarrow \PMod(\Sigma,\P_n) \longrightarrow \M(\Sigma,\P_n) \longrightarrow \S_n\longrightarrow 1,
\]
where the group $\S_n$ is generated by the images of half-twists in $\M(\Sigma,\P_n)$ permuting points of $\P_n$ (see e.g.~\cite{LP}). In the case of the torus with three punctures we can choose representatives of $\S_3$ as rigid rotations of the torus which send $\P_3$ to itself. We conclude that the above sequence splits: $\M(\Sigma_{1,0},\P_3)\cong \PMod(\Sigma_{1,0},\P_3)\rtimes\S_3$. (An explicit relation between the rotations $\sigma_1,\sigma_2$ and the corresponding half-twists which induce the same permutation of the punctures is given by formulas~\eqref{eq:tausigma} in Section~\ref{sec:6}.) Similarly, the orientation-reversing reflection forms an order $2$ subgroup in $\operatorname{Homeo}^{\pm}(\Sigma_{1,0},\P_3)$, hence $\M^*(\Sigma_{1,0},\P_3)$ splits as well: $\M^*(\Sigma_{1,0},\P_3)\cong \M(\Sigma_{1,0},\P_3)\rtimes \Z_2$. If we take as the generator of $\Z_2$ the reflection $\iota$ in the horizontal plane in Figure~\ref{fig:torus}, we observe that $\iota$ leaves every puncture fixed and commutes with all of $\S_3$, and hence we have the decomposition: $\M^*(\Sigma_{1,0},\P_3)\cong \PMod(\Sigma_{1,0},\P_3)\rtimes(\S_3\times \Z_2)$. The commutation action formulas given in the theorem are checked by direct inspection.

The fact that $\PMod(\Sigma_{1,0},\P_3)$ is isomorphic to $\ov{A[D_4]}$ with the given presentation on Dehn twist generators as the standard Artin generators for $\ov{A[D_4]}$ was shown in~\cite[Corollary\,9]{S}. 
It follows from \mbox{Korkmaz's} work~\cite{Kor} that $\Com(\M^*(\Sigma_{1,0},\P_3))$ is isomorphic to $\M^*(\Sigma_{1,0},\P_3)$, see~\cite[Theorem\,3.1]{BM}. Since $\ov{A[D_4]}$ is a subgroup of finite index in $\M^*(\Sigma_{1,0},\P_3)$, we have:
\[
\Com(\ov{A[D_4]})\cong \Com(\M^*(\Sigma_{1,0},\P_3))\cong \M^*(\Sigma_{1,0},\P_3).\qedhere
\]
\end{proof}

As a byproduct we get a description of the automorphism groups of $A[D_4]$ and $\ov{A[D_4]}$, which answers the question of Crisp and Paris for the case of $D_4$~\cite{CrispParisBD}.
\begin{cor}\label{cor:aut} \ 
\begin{enumerate} 
\item $\Aut(A[D_4])\cong\Aut(\ov{A[D_4]})\cong\ov{A[D_4]}\rtimes(\S_3\times\Z_2)\cong \M^*(\Sigma_{1,0},\P_3)$;
\item $\Out(A[D_4])\cong\Out(\ov{A[D_4]})\cong \S_3\times\Z_2$;
\item for $G=A[D_4]$, $\ov{A[D_4]}$ the following short exact sequence splits:
\[
1\longrightarrow\Inn(G)\longrightarrow\Aut(G)\longrightarrow\Out(G)\longrightarrow1.
\]
\end{enumerate}
\end{cor}
\begin{proof}
For part (1), we find $\Aut(\ov{A[D_4]})$ first. By Theorem~\ref{thm:comd4}, $\ov{A[D_4]}$ is a subgroup of finite index in $\M^*=\M^*(\Sigma_{1,0},\P_3)$. We notice that each automorphism $\alpha$ of $\ov{A[D_4]}$ can be treated as an injective homomorphism of a finite index subgroup ${}\ov{A[D_4]}$ of $\M^*$ into $\M^*$. By the Main Theorem 1.2 of~\cite{BM}, any injection of a finite index subgroup of $\M^*$ into $\M^*$ is induced by a conjugation with some element of $\M^*$. Since $\ov{A[D_4]}$ is normal in $\M^*$, the automorphisms of $\ov{A[D_4]}$ are induced by conjugation with the whole of $\M^*$. We are going to show that the conjugation action of $\M^*$ on ${}\ov{A[D_4]}$ is faithful. 

Let $\ell\colon A[D_4]\longrightarrow\Z$ be the homomorphism which maps each generator of $A[D_4]$ to $1$. Denote $d=\ell(\Delta)=12$, where $\Delta=(a_1a_2a_3a_4)^3$ is the generator of the center of $A[D_4]$ (see~\cite{BrieskornSaito}). Clearly, the center of $A[D_4]$ is being mapped onto $d\Z$ by $\ell$, hence we can define $\bar\ell\colon\ov{A[D_4]}\longrightarrow\Z/d\Z$ by reducing $\ell$ modulo $d$. Denote as $i_g$ the conjugation of $\M^*$ by an element $g\in\M^*$. We notice that if $g\in\ov{A[D_4]}\rtimes\S_3$, then $\bar\ell\circ i_g=\bar\ell$, whereas if $g=\iota$, then $\bar\ell\circ i_g=-\bar\ell$. Thus, if an element $g\in\M^*$ centralizes all of ${}\ov{A[D_4]}$, then $g$ actually belongs to $\M=\ov{A[D_4]}\rtimes\S_3$. Furthermore, such an element $g$ will centralize all Dehn twists about the curves $a_1,\dots,a_4$ from Figure~\ref{fig:torus}, and hence fix the isotopy classes of these curves up to orientation, by~\cite[Fact\,3.6]{FM}. Now notice that these curves cut the surface $(\Sigma_{1,0},\P_3)$ into punctured disks, and each such disk is uniquely determined by the subset of these curves that form its boundary. We conclude that $g$ fixes each one (of the isotopy classes rel $\P_3$) of these punctured disks setwise. Moreover $g$ either induces the trivial mapping class on each of these punctured disks, or an orientation reversing involution. Since $g$ was orientation preserving, we conclude that $g$ induces the trivial mapping class on each of these disks, and hence $g$ induces the identity automorphism of the graph formed by the representatives of the isotopy classes of oriented curves $a_1$, \dots, $a_4$. By the Alexander method (\cite[Prop.\,2.8]{FM}), $g$ must be trivial. This proves that the kernel of the conjugation action of $\M^*$ on $\ov{A[D_4]}$ is trivial, and hence that $\Aut(\ov{A[D_4]})\cong\M^*$.

By Proposition~8 of~\cite{CharneyCrisp}, we have $\Aut(A[D_4])\cong\Aut^*(\ov{A[D_4]})$, where $\Aut^*(\ov{A[D_4]})$ is the subgroup of all automorphisms $\alpha\in\Aut(\ov{A[D_4]})$ such that $\bar\ell\circ\alpha=\pm\bar\ell$. We established above that $\Aut^*(\ov{A[D_4]})$ is actually the whole group $\Aut(\ov{A[D_4]})$, and hence $\Aut(A[D_4])\cong\Aut(\ov{A[D_4]})$.

To prove part (2) we observe that $\Inn(A[D_4])\cong\ov{A[D_4]}\cong\Inn(\ov{A[D_4]})$ since the center of $\ov{A[D_4]}$ is trivial. Part (3) follows from parts (1) and (2).
\end{proof}

\section{Torsion in Artin groups of spherical type modulo their centers}\label{sec:4}
We will need a description of finite order elements in groups $\ov{A[D_4]}$, $\ov{A[F_4]}$, and $\ov{A[H_4]}$. Such description is probably known to experts, although it is not easily available in the literature. For groups of spherical types $A_n$, $B_n$, $D_n$, and $I_2(m)$, this information is present in~\cite[Sect.\,4.4]{LL} (and for type $A_n$ has been known long ago, see references therein). We complement it with the data on the rest of irreducible Artin groups of spherical type in Theorem~\ref{th:orders} below, which we distilled from several sources and from our own computations using the package \textsf{CHEVIE}~\cite{Chevie,Chevie0} (based on the computer algebra system \textsf{GAP3}~\cite{GAP3}), since we believe that this result may have an independent interest. 

It turns out that for each spherical type $\Gamma$ there exists a finite set of (one, two or three) basic torsion elements $\varepsilon_p$ such that all other torsion elements of $\ov{A[\Gamma]}$ are conjugate to a suitable power of one of $\varepsilon_p$. Namely, we have the following Theorem.

\begin{thm}\label{th:orders}
Let $A[\Gamma]$ be the irreducible Artin group of type $\Gamma$ (other than $A_1$), and $\ov{A[\Gamma]}$ be the quotient of $A[\Gamma]$ by its center. Then:
\begin{enumerate}
\item An element of finite order $d>1$ exists in $\ov{A[\Gamma]}$ if and only if $d$ is listed in the second column of Table~\ref{tab:1}.
\item Any torsion element in $\ov{A[\Gamma]}$ is conjugate to a power of one of the basic torsion elements $\varepsilon_p$ listed in the third column of Table~\ref{tab:1}. 
\item For any two basic elements $\varepsilon_p$ and $\varepsilon_q$ of orders $p$ and $q$, respectively, from Table~\ref{tab:1}, if $d=\gcd(p,q)$, then the power $\varepsilon_p^{p/d}$ is conjugate to $\varepsilon_q^{q/d}$ in $\ov{A[\Gamma]}$. Powers of the same element $\varepsilon_p^s$ and $\varepsilon_p^t$ are not conjugate in $\ov{A[\Gamma]}$, if $s\ne t\pmod p$.
\end{enumerate}
\end{thm}

\begin{table}[hbt!]
\begin{center}
{\renewcommand{\arraystretch}{0}%
\begin{tabular}{|c|c|c|c|}
\hline\strut
type $\Gamma$ & orders of torsion & basic torsion elements & relations\\
\hline \rule{0pt}{1pt}&\\
\hline\strut $A_n$, $n\ge2$ & all divisors of $n$, $n+1$ & 
			\parbox{5.5cm}{\centering\strut$\varepsilon_{n+1}=s_1s_2\dots s_n$\\
									\strut$\varepsilon_{n}=s_1(s_1s_2\dots s_n)$}&\\[3mm]
\hline\strut $B_n$, $n\ge2$ & all divisors of $n$ & \strut$\varepsilon_{n}=s_1s_2\dots s_n$&\\
\hline\strut $D_n$, $n$ even $\ge4$ & all divisors of $n-1$, $n/2$ & 
			\parbox{5.5cm}{\centering\strut$\varepsilon_{n-1}=s_ns_{n-1}\dots s_3s_2s_1$\\
										\strut$\varepsilon_{n/2}=(s_n\dots s_3)s_2(s_n\dots s_3)s_1$}&\\[3mm]
\hline\strut $D_n$, $n$ odd $\ge5$ & all divisors of $2n-2$, $n$ & 
			\parbox{5.5cm}{\centering\strut$\varepsilon_{2n-2}=s_ns_{n-1}\dots s_3s_2s_1$\\
										\strut$\varepsilon_{n}=(s_n\dots s_3)s_2(s_n\dots s_3)s_1$}&\\[3mm]
\hline\strut $E_6$ & 2,\,3,\,4,\,6,\,8,\,9,\,12 &	
\parbox{5.5cm}{\centering\strut$\varepsilon_{12}=
																									s_4s_2\cdot s_3s_1\cdot s_5s_6$\\
								\strut$\varepsilon_{9}=
																				s_4s_2\cdot s_5s_4s_3s_1\cdot s_6s_5$\\
								\strut$\varepsilon_{8}=
																				s_4s_3s_1\cdot s_5s_4s_2\cdot s_3s_6s_5$}& 
\parbox{2cm}{\centering\strut$\varepsilon_{12}^4=\varepsilon_{9}^3$\\[1mm]
						\strut$\varepsilon_{12}^3=\varepsilon_{8}^2$}\\[3mm]				
\hline\strut $E_7$ & 3,\,7,\,9 &	
\parbox{5.5cm}{\centering\strut$\varepsilon_{9}=
																								s_4s_2\cdot s_3s_1\cdot s_5s_6s_7$\\
								\strut$\varepsilon_{7}=
																				s_4s_2\cdot s_7s_6s_5\cdot s_4s_2\cdot s_3s_1$}&\\[3mm]
\hline\strut $E_8$ & 2,\,3,\,4,\,5,\,6,\,10,\,12,\,15 &	
\parbox{5.8cm}{\centering\strut$\varepsilon_{15}=
																									s_4s_2\cdot s_3s_1\cdot s_8s_7s_6s_5$\\
								\strut$\varepsilon_{12}=
																				s_4s_2\cdot s_3s_1\cdot s_4s_3\cdot s_8s_7s_6s_5$\\
								\strut$\varepsilon_{10}=
																				s_4s_2\cdot s_3s_1\cdot s_6s_5s_4s_3\cdot s_8s_7s_6s_5$}& \\[3mm]
\hline\strut $F_4$ & 2,\,3,\,4,\,6 &	
\parbox{5.5cm}{\centering\strut$\varepsilon_{6}=s_1s_2s_3s_4$\\
								\strut$\varepsilon_{4}=s_1s_2s_3s_4\cdot s_2s_3$}&
									$\varepsilon_6^3=\varepsilon_4^2$
								\\[3mm]				
\hline\strut $H_3$ & 3,\,5 &	
\parbox{5.5cm}{\centering\strut$\varepsilon_{5}=s_1s_2s_3$\\
								\strut$\varepsilon_{3}=s_1s_2\cdot s_1s_2s_3$}&\\				
\hline\strut $H_4$ & 2,\,3,\,5,\,6,\,10,\,15 &	
\parbox{5.5cm}{\centering\strut$\varepsilon_{15}=s_1s_2s_3s_4$\\
								\strut$\varepsilon_{10}=s_1s_2\cdot s_1s_2s_3s_4$\\
								\strut$\varepsilon_{6}=s_1s_2\cdot s_1s_2s_3s_2\cdot s_1s_2s_3s_4$}& \\[3mm]
\hline\strut $I_2(m)$, $m$ even $\ge6$ & all divisors of $m/2$ &	$\varepsilon_{m/2}=s_1s_2$&\\
\hline\strut $I_2(m)$, $m$ odd $\ge5$ & all divisors of $2$, $m$ &	
\parbox{5.5cm}{\centering\strut$\varepsilon_{m}=s_1s_2$\\
								\strut$\varepsilon_{2}=s_1(s_2s_1)^{(m-1)/2}$}&\\[3mm]
\hline
\end{tabular}}\end{center}
\caption{Orders of torsion and representatives of basic torsion elements up to conjugacy in groups~$\protect\ov{A[\Gamma]}$. An element $\varepsilon_p$ has order $p$, and is written in the standard generators $s_i$ of $\protect\ov{A[\Gamma]}$ using the numbering presented in Figure~\ref{fig:sph}.\label{tab:1}}
\end{table}

\begin{rem}
The basic elements $\varepsilon_p$ from Table~\ref{tab:1} were chosen as positive words of minimal length in the standard generators for $\ov{A[\Gamma]}$. In some cases the expressions for the generators $\varepsilon_p$ and $\varepsilon_q$ could be chosen with the extra compatibility condition that $\varepsilon_p^{p/d}$ is actually equal to $\varepsilon_q^{q/d}$ for $d=\gcd(p,q)$. When this was possible, we listed this information in the fourth column. Interestingly, for groups of types $E_8$ and $H_4$ such choices were not possible, as we checked using \textsf{CHEVIE} by computing all positive conjugates of given elements in the corresponding Artin groups $A[\Gamma]$.
\end{rem}

\begin{proof}[Proof of Theorem~\ref{th:orders}]
We will prove several claims, from which the Theorem will follow. Recall that an element of $A[\Gamma]$ is called \textit{periodic} if some power of it lies in the center. We will call a periodic element of $A[\Gamma]$ \textit{primitive}, if it is not a nontrivial power of any other element of $A[\Gamma]$. Let $Z$ denote the infinite cyclic center of $A[\Gamma]$, and let $\delta\in Z$ denote its standard generator.

\textit{Claim 1. Every element of finite order $\bar u\in\ov{A[\Gamma]}$ has a representative $u\in A[\Gamma]$ which is a power of some primitive periodic element $v\in A[\Gamma]$.}
Indeed, it was proved in~\cite[Prop.\,3.4]{Paris1} that there is an upper bound on finite orders of elements in any given $\ov{A[\Gamma]}$ of spherical type. This means that for an arbitrary representative $u\in A[\Gamma]$ such that $\bar u=uZ$, one cannot take roots of arbitrarily large degrees out of $u$ in $A[\Gamma]$, hence there exists a root $v$ of $u$ of maximal degree, which will be primitive.

\textit{Claim 2. All primitive periodic elements $v\in A[\Gamma]$ are roots of the standard generator $\delta$ of the center $Z$.} 
This was proved by Lee and Lee in~\cite[Th.\,3.14]{LL}.

Recall that the standard generator $\delta$ of the center is equal to $\Delta^\kappa$, where $\Delta$ is the Garside element of $A[\Gamma]$, with $\kappa=2$ for $\Gamma$ of the following types:
\begin{equation*}\label{Artin:inv}
A_n\,\, (n\ge2),\quad D_n\,\, \text{($n$ odd)},\quad E_6,\quad I_2(m)\,\, \text{($m$ odd)},\tag{*}
\end{equation*}
and $\kappa=1$ otherwise (see~\cite{BrieskornSaito}). 

We are going to get a description of which $d$-th roots of $\delta$ exist in $A[\Gamma]$ and to prove that they are all conjugate in $A[\Gamma]$. These results can be distilled from works of Bessis and Springer. In~\cite[Th.\,12.4]{Bes} Bessis proved for a certain central element $\iota=\delta_0^h\in A[\Gamma]$, where $\delta_0$ the Garside element of the \emph{dual} Garside structure on the Artin group of spherical type $A[\Gamma]$, and $h$ the Coxeter number of the Coxeter group $W[\Gamma]$, the following:

\textit{Claim 3.
For $d>1$, a $d$-th root of $\iota$ exists in $A[\Gamma]$ if and only if $d$ is a \emph{regular number} in the sense of Springer~\cite{Sp}. 
If a $d$-th root of $\iota$ exists, it is unique up to conjugation in $A[\Gamma]$.}

The exact relationship between $\delta$ and $\iota$ is captured by Theorem~12.3 of~\cite{Bes}, which states that $\delta$ is equal to $\delta_0^{h'}$ where $h'=h/\gcd(d_1,\dots,d_n)$, and $d_1,\dots,d_n$ are the degrees of the Coxeter group $W[\Gamma]$. By direct inspection of the lists of degrees for all spherical Coxeter groups~\cite[p.\,59]{Humph1} we observe that $\gcd(d_1,\dots,d_n)=1$ exactly for Artin groups of types~\eqref{Artin:inv} above, and it equals $2$ otherwise. Since $\iota=\delta_0^h$, we conclude that: 

\textit{Claim 4. The element $\iota$ equals $\Delta^2$ in all cases.}

Now we are ready to prove the Theorem. For part (1), we point out that in~\cite[Sect.\,5]{Sp} Springer listed explicitly the regular numbers for all Artin groups of spherical type (except $I_2(m)$, for which this information can be taken from~\cite[Table~1]{LL}). Combined with the first part of Claim~3, this gives us the description of all roots of $\Delta^2$ up to conjugation in $A[\Gamma]$. In the cases when $\delta=\Delta^2$, i.e.\ for groups in the list~\eqref{Artin:inv} above, this gives the final answer: the orders of torsion elements of $\ov{A[\Gamma]}$ are exactly the regular numbers from~\cite[Section\,5]{Sp} and \cite[Table\,1]{LL}. 

However for groups not in the list~\eqref{Artin:inv}, we need to adjust the lists of regular numbers so that they correspond to roots of $\delta=\Delta$, and not of $\Delta^2$ (cf.~\cite{Sh}). For that, select the subset of regular numbers which are maximal with respect to divisibility, say $r_1,\dots,r_k$. These numbers correspond to primitive periodic elements. We notice that all numbers $r_1,\dots,r_k$ will be divisible by $2$ since by Claim~2 all primitive periodic elements should be divisors of $\delta$. By direct inspection of the data in~\cite[Section\,5]{Sp} and \cite[Table\,1]{LL} we notice that all numbers $r_i$ are bigger than $2$. Form a set $r_1/2,\dots,r_k/2$. These numbers will be the orders of basic torsion elements in $\ov{A[\Gamma]}$, and the orders of all other torsion elements will be all the divisors of these. This completes the description of orders of torsion elements in groups ${}\ov{A[\Gamma]}$, which we list in the second column of Table~\ref{tab:1}. 

Now we prove part (2). Let $\bar u$ be an element of finite order in $\ov{A[\Gamma]}$. By Claim~1, $\bar u$ has a representative $u\in A[\Gamma]$ such that $u=v^k$ for some primitive periodic element $v\in A[\Gamma]$ and some $k\in\Z$. By Claim~2, $v$ is a root of $\delta$, and hence of $\iota$. Let $\epsilon_p$ denote the element of $A[\Gamma]$ corresponding to an element $\varepsilon_p\in\ov{A[\Gamma]}$ from the third column of Table~\ref{tab:1}, i.e.\ $\epsilon_p$ is given by the same word in the standard generators of $A[\Gamma]$ as $\varepsilon_p$. We claim that each $\epsilon_p$ is a root of the $p$-th degree of $\delta$. This was proved for the infinite series of types $A_n$, $B_n$, $D_n$, $I_2(m)$ in~\cite[Lemma\,4.3]{LL}, and for standalone systems of types $E_6$, $E_7$, $E_8$, $F_4$, $H_3$, $H_4$ it can be checked using \textsf{CHEVIE}. We observe that, for a given $A[\Gamma]$, the degrees $p$ of elements $\epsilon_p$ form the set of all maximal elements with respect to divisibility of the numbers listed in the second column of Table~\ref{tab:1}. Thus we conclude that elements $\epsilon_p$ are primitive roots of $\delta$ (and hence of $\iota$), and account for all possible degrees of primitive roots of $\delta$. Thus, by the second part of Claim~3, $v$ is conjugate to one of the elements $\epsilon_p$, and hence $\bar u=v^kZ$ is conjugate to the $k$-th power of the corresponding element $\varepsilon_p$. This proves (2).

Similarly, if $d=\gcd(p,q)$, then both $\epsilon_p^{p/d}$ and $\epsilon_q^{q/d}$ are roots of $\delta$ of the same degree $d$, and hence they are roots of $\iota$ of equal degree. By the second part of Claim~3, these elements are conjugate in $A[\Gamma]$ and hence $\varepsilon_p^{p/d}$ and $\varepsilon_q^{q/d}$ are conjugate in $\ov{A[\Gamma]}$, which proves the first sentence of part~(3). To prove the second sentence of part (3), consider the homomorphism $\ell\colon A[\Gamma]\longrightarrow \Z$, which sends every standard generator to $1$, and its reduction modulo $\ell(\delta)$, which gives rise to $\bar\ell\colon\ov{A[\Gamma]}\longrightarrow\Z/\ell(\delta)\Z$. Clearly, conjugation in $A[\Gamma]$ and $\ov{A[\Gamma]}$ leaves $\ell$ and $\bar\ell$ invariant. Now to prove that if $s\ne t\pmod p$, the elements $\varepsilon_p^s$ and $\varepsilon_p^t$ are not conjugate, we may assume without loss of generality that $0\le s,t\le p-1$ and consider the value of $\bar\ell$ on these elements. We observe from the explicit expressions for elements $\epsilon_p$ that $\ell(\delta)=p\ell(\epsilon_p)$ in each case, and hence we have: $\ell(\epsilon_p^s)=s\ell(\epsilon_p)<\ell(\delta)$ and $\ell(\epsilon_p^t)=t\ell(\epsilon_p)<\ell(\delta)$. Therefore $\bar\ell(\varepsilon_p^s)\ne \bar\ell(\varepsilon_p^t)$, and part (3) is proved.
\end{proof}
We can now describe conjugacy classes of torsion elements in groups $\ov{A[\Gamma]}$. In the following Corollary, $\phi$ denotes Euler's totient function.
\begin{cor}\label{cor:phi}
Let $\ov{A[\Gamma]}$ be as in Theorem~\ref{th:orders}. Then torsion elements of order $d>1$ in $\ov{A[\Gamma]}$, if they exist, form $\phi(d)$ conjugacy classes, with representatives of the form $\big(\varepsilon_p^{p/d}\big)^\ell$, where $\varepsilon_p$ is any basic element from the third column of Table~\ref{tab:1} such that $d$ divides $p$, and $\ell$ runs through all positive integers less than $d$ which are coprime to $d$.
\end{cor}
\begin{proof}
By part (2) of Theorem~\ref{th:orders}, any torsion element $x$ in $\ov{A[\Gamma]}$ is conjugate to an element of the form $g=\varepsilon_p^k$, where $\varepsilon_p$ is a basic torsion element from the third column of Table~\ref{tab:1}. Since $\varepsilon_p$ has order $p$, we conclude that the order of $\varepsilon_p^k$ is equal to $p/\gcd(p,k)$. In particular, if $d$ is the order of a torsion element in $\ov{A[\Gamma]}$, $d$ divides $p$. If $g$ has order $d$, we have: $g^d=\big(\varepsilon_p^k\big)^d=1$, so that $kd=\ell p$ for some $\ell\ge 1$. Since $d$ divides $p$, we conclude that $g=\varepsilon_p^k=\big(\varepsilon_p^{p/d}\big)^\ell$. Let $c=\gcd(\ell,d)$. If $c\ne1$ then $g^{d/c}=\big(\varepsilon_p^k\big)^{d/c}=\big(\varepsilon_p^{p/d}\big)^{\ell d/c}=\big(\varepsilon_p^{p}\big)^{\ell/c}=1$ and $d$ is not the order of $g$. Hence, $\ell$ and $d$ are coprime. Notice also that $\big(\varepsilon_p^{p/d}\big)^\ell=\big(\varepsilon_p^{p/d}\big)^{\ell+md}$ for any $m\in\Z$, so that we may assume that $\ell$ belongs to $\{1,\dots,d-1\}$. If $\ell_1\ne\ell_2$ are two elements from $\{1,\dots,d-1\}$ coprime to $d$ then we claim that $\frac{p}{d}\cdot\ell_1\ne \frac{p}{d}\cdot \ell_2\pmod p$. Indeed, if $\frac{p}{d} (\ell_1-\ell_2)=0\pmod p$ this means that $\ell_1-\ell_2=md$ for some $m\in\Z$. But since $0<\ell_1,\ell_2<d$, we must have $m=0$ and $\ell_1=\ell_2$. This contradiction shows that $\frac{p}{d}\cdot\ell_1\ne \frac{p}{d}\cdot \ell_2\pmod p$, and by part~(3) of Theorem~\ref{th:orders}, all elements $\big(\varepsilon_p^{p/d}\big)^\ell$ are mutually non-conjugate when $\ell$ runs through all positive integers less than $d$ and coprime to $d$.

If there exists another basic torsion element $\varepsilon_q$ such that the torsion element $x$ is conjugate to some power $h=\varepsilon_q^m$, then the above analysis shows that $h=\big(\varepsilon_q^{q/d}\big)^\ell$. By part~(3) of Theorem~\ref{th:orders}, $\varepsilon_q^{q/d}$ is conjugate to $\varepsilon_p^{p/d}$. Thus, any choice of basic element $\varepsilon_p$ such that $d$ divides $p$ will work for description of representatives for conjugacy classes of elements of order $d$.
\end{proof}

\section{Non-commensurability of Artin groups of types \texorpdfstring{$H_4$}{H4} and \texorpdfstring{$D_4$}{D4}}\label{sec:5}
\begin{prop}\label{prop:d4h4}
The groups $A[H_4]$ and $A[D_4]$ are not commensurable.
\end{prop}
\begin{proof}
Suppose that $A[H_4]$ is commensurable with $A[D_4]$. Then by Theorem~\ref{thm:comm}, $\ov{A[H_4]}$ and ${}\ov{A[D_4]}$ are also commensurable, and moreover $\ov{A[H_4]}$ embeds into $\Com(\ov{A[H_4]})\cong\Com(\ov{A[D_4]})$, which is isomorphic to $\ov{A[D_4]}\rtimes (\S_3\times \Z_2)$ due to Theorem~\ref{thm:comd4}. Denote $\Phi\colon \ov{A[H_4]}\lhook\joinrel\longrightarrow \ov{A[D_4]}\rtimes (\S_3\times \Z_2)$ such an embedding, and consider the natural projection $\theta\colon \ov{A[D_4]}\rtimes (\S_3\times \Z_2) \longrightarrow \S_3\times \Z_2$. It is known~\cite{BrieskornSaito} that the center of $A[H_4]$ is generated by the element $\Delta=(s_1s_2s_3s_4)^{15}$, where $s_i$'s are the standard generators of $A[H_4]$ numbered as in Figure~\ref{fig:sph}. By considering the homomorphism $\ell\colon A[H_4]\longrightarrow\Z$ sending each standard generator to $1$, we see that the center of $A[H_4]$ maps onto a subgroup $60\Z$, and hence no power less than $15$ of the element $\mu=s_1s_2s_3s_4$ lies in the center of $A[H_4]$. This shows that the image $\bar\mu$ of the element $\mu$ in $\ov{A[H_4]}$ has order exactly $15$. Let $M=\langle \bar\mu\rangle$ be the cyclic subgroup of order $15$ in $\ov{A[H_4]}$ generated by $\bar\mu$. Clearly, $M\cong\Z_3\oplus \Z_5$. The subgroup $\Phi(\Z_5)$ of $\Phi(M)$ lies in the kernel of $\theta$, since $\Phi$ is injective and $\S_3\times \Z_2$ has no elements of order $5$. But $\Ker\theta$ lies in $\ov{A[D_4]}$, which does not have elements of order $5$ either, since the maximal order of torsion elements in $\ov{A[D_4]}$ is $3$, see Table~\ref{tab:1}. This contradicts the injectivity of $\Phi$ and proves that the Artin groups $A[H_4]$ and $A[D_4]$ are not commensurable.
\end{proof}

The following Corollary gives partial answer to Question~\ref{que:emb} concerning which central quotients ${}\ov{A[\Gamma]}$ admit embeddings into mapping class groups of punctured surfaces.

\begin{cor}\label{cor:h4mcg}
The group $\ov{A[H_4]}$ does not admit any embedding into the extended mapping class group $\M^*(\Sigma_{1,0},\P_3)$.
\end{cor}
\begin{proof}
Indeed, in the proof of Proposition~\ref{prop:d4h4} we showed that the assumption that there exists an embedding $\Phi\colon \ov{A[H_4]}\lhook\joinrel\longrightarrow {}\ov{A[D_4]}\rtimes (\S_3\times \Z_2)$ leads to a contradiction. Since by Theorem~\ref{thm:comd4}, \mbox{$\ov{A[D_4]}\rtimes (\S_3\times \Z_2)$} is isomorphic to $\M^*(\Sigma_{1,0},\P_3)$, the Corollary follows.
\end{proof}

\section{Non-commensurability of Artin groups of types \texorpdfstring{$F_4$}{F4} and \texorpdfstring{$D_4$}{D4}}\label{sec:6}
\begin{prop}\label{prop:d4f4}
The groups $A[F_4]$ and $A[D_4]$ are not commensurable.
\end{prop}
\begin{proof}
Suppose that $A[F_4]$ is commensurable with $A[D_4]$. Then by Theorem~\ref{thm:comm}, $\ov{A[F_4]}$ and $\ov{A[D_4]}$ are commensurable, and $\ov{A[F_4]}$ embeds into $\Com(\ov{A[F_4]})\cong\Com(\ov{A[D_4]})$, isomorphic to $\M^*=\M^*(\Sigma_{1,0},\P_3)\cong\ov{A[D_4]}\rtimes (\S_3\times \Z_2)$ by Theorem~\ref{thm:comd4}.  Let $\Phi\colon \ov{A[F_4]}\lhook\joinrel\longrightarrow\M^*$ be such an embedding, and denote $K=\ov{\CA[D_4]}$. We are going to classify all projections of $\Phi$ to $\M^*/K$ and obtain a contradiction with the injectivity of $\Phi$.

In what follows the standard generators of $\ov{A[D_4]}$ will be denoted as $a_1,a_2,a_3,a_4$ as in Figure~\ref{fig:torus}, and the standard generators of $\ov{A[F_4]}$ as $s_1,s_2,s_3,s_4$, with the numeration as in Figure~\ref{fig:sph}. 

Recall that $K=\ov{\CA[D_4]}$ is normally generated in $\ov{A[D_4]}$ by the squares of the standard generators $\{a_1^2, a_2^2, a_3^2, a_4^2\}$. The conjugation by $\S_3$ leaves this set invariant, whereas the conjugation by $\Z_2=\langle\iota\rangle$ sends all these elements to their inverses. This implies that $K$ is normal in $\M^*$ and we can consider the natural projections 
\[
\Theta\colon\M^*\longrightarrow \M^*/K,\text{\quad and \quad} \varphi=\Theta\circ\Phi.
\]
We notice that $\M^*/K=(\ov{A[D_4]}/K)\rtimes (\S_3\times\Z_2)=\ov{W[D_4]}\rtimes(\S_3\times\Z_2)$, where $\ov{W[D_4]}$ is the quotient of the Coxeter group $W[D_4]$ by its center. Indeed, the image of $\Theta$ has a presentation with the generators $a_1,\dots,a_4$ and the relations listed in Theorem~\ref{thm:comd4}, with the additional relations $a_1^2=a_2^2=a_3^2=a_4^2=1$. According to~\cite[Corollary\,3.19]{Humph1}, the center of $W[D_4]$ is generated by the image of $\Delta=(a_1a_2a_3a_4)^3$, the standard generator of the center of $A[D_4]$, so modding out $\Delta$ and the squares of the generators $a_i^2$ can be performed in any order. Notice also that the generator of $\Z_2$ acts trivially on $\ov{W[D_4]}$, so we actually have: $\M^*/K=\big(\ov{W[D_4]}\rtimes \S_3\big)\times\Z_2$, and hence $\varphi$ decomposes as $\varphi=(\psi,\zeta)$ for some $\psi\colon\ov{A[F_4]}\longrightarrow\ov{W[D_4]}\rtimes\S_3$ and $\zeta\colon\ov{A[F_4]}\longrightarrow\Z_2$.

We can describe all possible $\zeta$ easily. Identify $\Z_2$ with $\{1,\iota\}$. Since $s_1s_2s_1=s_2s_1s_2$ in $\ov{A[F_4]}$, we must have $\zeta(s_1)=\zeta(s_2)$, and similarly, $\zeta(s_3)=\zeta(s_4)$. These two conditions define four possible homomorphisms $\zeta_i$ given on the standard generators as follows: 
\begin{align*}
\zeta_1\colon(s_1,s_2,s_3,s_4)&\longmapsto (1,1,1,1); & \zeta_3\colon(s_1,s_2,s_3,s_4)&\longmapsto (1,1,\iota,\iota); \\
\zeta_2\colon(s_1,s_2,s_3,s_4)&\longmapsto (\iota,\iota,1,1); & \zeta_4\colon(s_1,s_2,s_3,s_4)&\longmapsto (\iota,\iota,\iota,\iota).
\end{align*}

Taking into account torsion elements in $\ov{A[F_4]}$ imposes a significant restriction on possible homomorphisms $\varphi$. Let $v=s_1s_2s_3s_4\in\ov{A[F_4]}$. By reasoning as in the proof of Proposition~\ref{prop:d4h4}, we conclude that $v$ has order $6$ in $\ov{A[F_4]}$. By Proposition~\ref{prop:ca}, $K$ embeds into $\CA[D_4]$, which is a torsion free group (since it is bi-orderable). Hence the subgroup $\langle v\rangle$ maps isomorphically into $\M^*/K$, i.e.\ $\varphi(v)$ must have order~$6$. Notice that all homomorphisms $\zeta_i$ map $v$ to $1$, which means that $\langle v\rangle$ maps isomorphically into $\ov{W[D_4]}\rtimes\S_3$ and the order of $\psi(v)$ is~$6$ as well.

To describe all possible homomorphisms $\psi\colon\ov{A[F_4]}\longrightarrow\ov{W[D_4]}\rtimes\S_3$ we use the computer algebra system \textsc{Magma}~\cite{Magma}, which has a functionality of computing all homomorphisms, including non-surjective ones, from a finitely presented group to a finite group, up to the conjugacy in the latter (via the command \texttt{Homomorphisms} with the flag \texttt{Surjective:=false}). Performing such a computation we get that there exist 286 conjugacy classes of possible homomorphisms $\psi$. Out of these classes, only 10 satisfy the condition that $\psi(v)$ has order~$6$. This number is further reduced to five if we consider these homomorphisms up to the graph automorphism of $\ov{A[F_4]}$ interchanging $s_1\longleftrightarrow s_4$, $s_2\longleftrightarrow s_3$. 

Analyzing these five cases, we observe that in four of them the following property holds: $\psi(s_1)=\psi(s_2)$ and the order of this element is either $1$ or $2$. Hence the same will be true for $\varphi=(\psi,\zeta)$ where $\zeta$ is any of the four homomorphisms $\zeta_i$ above, i.e.\ $\varphi(s_1)=\varphi(s_2)$ has order either $1$ or $2$. Now we can produce generalized torsion elements in $K$, and this will yield the non-faithfulness of $\Phi$. Indeed, take $\alpha=s_1s_2^{-1}$, $\beta=s_1s_2\in \ov{A[F_4]}$. From the above, we have: $\varphi(\alpha)=\varphi(\beta)=1$, hence $\alpha,\beta\in\Ker\Theta=K$. In $\ov{A[F_4]}$ (and even in $A[F_4]$) the following identity holds: $\alpha\cdot\beta\alpha\beta^{-1}\cdot\beta^2\alpha\beta^{-2}=1$, which can be checked either by hand or using the package \textsf{CHEVIE}. This means that $\Phi(\alpha)$ is a generalized torsion element in $K$, which is impossible due to the fact that $K$ lies in $\CA[D_4]$ and the group $\CA[D_4]$ is bi-orderable, by Proposition~\ref{prop:ca} and~\cite{RZ}.

In the remaining case out of the five classes of homomorphisms $\psi\colon\ov{A[F_4]}\longrightarrow\ov{W[D_4]}\rtimes\S_3$ found above, a representative can be chosen in the following form:
\[
\psi\colon(s_1,s_2,s_3,s_4)\longmapsto(\hat a_3,\hat a_2,\sigma_1,\sigma_2),
\]
where $\hat a_1,\dots,\hat a_4$ denote the images of the standard generators of $\ov{A[D_4]}$ in $\ov{W[D_4]}$. 
Observe that for each $i$ the homomorphism $\varphi=(\psi,\zeta_i)$ sends the squares of the standard generators of $\ov{A[F_4]}$ to~$1$, and hence $\Phi\big(\ov{\CA[F_4]}\big)\subset\Ker\Theta=K$. An easy check in \textsf{GAP}~\cite{GAP4} shows that, for each $i$, the order of the image of $\varphi=(\psi,\zeta_i)$ is $576$, which is exactly the order of $\ov{W[F_4]}=\ov{A[F_4]}/\ov{\CA[F_4]}$. This means that no element of $\ov{A[F_4]}\setminus\ov{\CA[F_4]}$ maps into $K$ via $\Phi$, and hence no generalized torsion can be detected in $\im\Phi\cap K$. Thus, to obtain a contradiction we need to proceed in some other way. 

We claim that the image $\Phi(\ov{A[F_4]})$ has finite index in $\M^*$. Indeed, by assumption, $\ov{A[F_4]}$ and ${}\ov{A[D_4]}$ are commensurable, hence some finite index subgroup $H_1\subseteq\ov{A[F_4]}$ is isomorphic to a finite index subgroup $H_2\subseteq\ov{A[D_4]}$. Since $\ov{A[D_4]}$ has finite index in $\M^*$, it follows that $\Phi(H_1)\subseteq\M^*$ is isomorphic to a finite index subgroup $H_2\subseteq\M^*$. By the Main Theorem~1.2 of~\cite{BM}, any injection of a finite index subgroup of $\M^*$ to $\M^*$ is induced by an inner automorphism of $\M^*$. This means that $\Phi(H_1)=gH_2g^{-1}$ for some $g\in\M^*$, and hence $\Phi(H_1)$ has finite index in $\M^*$. Thus the index of $\Phi(\ov{A[F_4]})$ in $\M^*$ is also finite.

Consider the automorphism $\alpha\colon\ov{A[F_4]}\longrightarrow\ov{A[F_4]}$ induced by the graph automorphism which interchanges $s_1\longleftrightarrow s_4$, $s_2\longleftrightarrow s_3$. Let $\Phi'=\Phi\circ\alpha\colon\ov{A[F_4]}\lhook\joinrel\longrightarrow\M^*$ be the corresponding embedding. Since $\Phi(\ov{A[F_4]})$ has finite index in $\M^*$, the same is true for $\Phi'(\ov{A[F_4]})$, and applying the Main Theorem~1.2 of~\cite{BM} again, we conclude that $\Phi'(x)=h\Phi(x)h^{-1}$ for some $h\in\M^*$ and all $x\in\ov{A[F_4]}$.

Observe now that since $\iota\sigma_i\iota=\sigma_i$ and $\iota a_i\iota=a_i^{-1}$, we conclude that for any $j=1,\dots,4$, $\varphi=(\psi,\zeta_j)$ has the following properties:
\[
\varphi(s_2s_1)=\hat a_2\hat a_3\text{\quad and\quad}\varphi(s_3s_4)=\sigma_1\sigma_2.
\]
This means that 
\[
\Phi(s_2s_1)\in a_2a_3K\subseteq \ov{A[D_4]}\text{\quad and\quad} \Phi(s_3s_4)\in \sigma_1\sigma_2 K \subseteq \M^*\setminus\ov{A[D_4]},
\]
since $\sigma_1\sigma_2$ is a nontrivial element of $\S_3$. But then
\[
\Phi'(s_2s_1)=\Phi(s_3s_4)\notin{}\ov{A[D_4]}\text{\quad and \quad}\Phi'(s_3s_4)=\Phi(s_2s_1)\in\ov{A[D_4]}, 
\]
which contradicts the fact that $\Phi'(x)=h\Phi(x) h^{-1}$, since $\ov{A[D_4]}$ is normal in $\M^*$.

This contradiction finishes the analysis of the last remaining case of $\psi$ and shows that an embedding $\Phi\colon\ov{A[F_4]}\lhook\joinrel\longrightarrow\M^*$, such that $\Theta\circ\Phi=(\psi,\zeta_i)$ and $\Phi(\ov{A[F_4]})$ has finite index in $\M^*$, does not exist. This implies that the groups $\ov{A[F_4]}$ and $\ov{A[D_4]}$ are not commensurable and hence that the groups $A[F_4]$ and $A[D_4]$ are not commensurable either.
\end{proof}


Remarkably, homomorphisms $\Psi\colon\ov{A[F_4]}\to \ov{A[D_4]}\rtimes\S_3$ such that their projection to $\ov{W[D_4]}\rtimes\S_3$ is given by $\psi\colon(s_1,s_2,s_3,s_4)\longmapsto(\hat a_3,\hat a_2,\sigma_1,\sigma_2)$, do exist. One such homomorphism is constructed in Example~\ref{ex:13} below. To describe it we need to recall the definition and basic properties of half-twists.

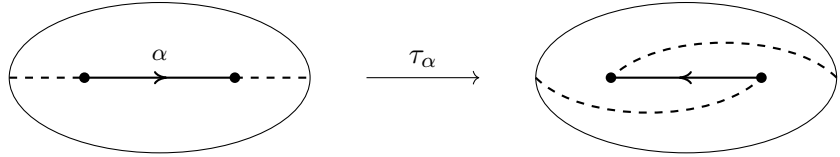
\begin{figure}[htb!]
\begin{center}
\begin{tikzpicture}
\draw (0,0) ellipse [x radius=2, y radius=1];
\fill (-1,0) circle (2pt);
\fill (1,0) circle (2pt);
\draw[dashed, thick] (-2,0)--(-1,0);
\draw[dashed, thick] (1,0)--(2,0);
\draw[->-=0.55,thick] (-1,0)--(1,0);
\draw (0,0.3) node {\footnotesize$\alpha$};

\begin{scope}[xshift=7cm]
\draw (0,0) ellipse [x radius=2, y radius=1];
\fill (-1,0) circle (2pt);
\fill (1,0) circle (2pt);
\draw[dashed, thick] (2,0) to [out=135,in=45,looseness=0.75] (-1,0); 
\draw[dashed, thick] (-2,0) to [out=-45,in=-135,looseness=0.75] (1,0);
\draw[->-=0.55,thick] (1,0)--(-1,0);
\end{scope}

\draw[->] (2.75,0) to (4.25	,0); 
\draw (3.5,0.25) node {$\tau_\alpha$};

\end{tikzpicture}
\end{center}
\caption{A right half-twist.\label{fig:halftwist}}
\end{figure}

Let $D$ be a twice punctured disk and $\alpha$ an embedded arc connecting the two punctures. The \emph{right half-twist along $\alpha$} is a homotopy class of homeomorphisms of $D$, which we denote $\tau_\alpha$, described by Figure~\ref{fig:halftwist}. If $D$ is embedded in a surface, we extend the right half-twist along $\alpha$ by the identity outside $D$ to obtain a mapping class of the surface, which we denote by $\tau_\alpha$ as well. If $\beta$ is an embedded arc between some other two punctures in the surface, then $\tau_\beta$ commutes with $\tau_\alpha$ if and only if $\alpha\cap \beta=\varnothing$, and $\tau_\beta\tau_\alpha\tau_\beta=\tau_\alpha\tau_\beta\tau_\alpha$ if and only if $\alpha$ and $\beta$ intersect in a common endpoint. Let $\gamma$ be a simple closed curve in the surface transversely intersecting the arc $\alpha$ once in its interior point, and let $a_\gamma$ be the \emph{right} Dehn twist about $\gamma$. Labru\`ere and Paris observed~\cite[Lemma\,2.3]{LP} that the Dehn twist $a_\gamma$ and the half-twist $\tau_\alpha$ satisfy the Artin relation of length four: 
\[
a_\gamma\tau_\alpha a_\gamma\tau_\alpha = \tau_\alpha a_\gamma\tau_\alpha a_\gamma.
\]

\begin{exa}\label{ex:13}
The above properties of half-twists allow us to engineer the following `geometric' homomorphism from $\ov{A[F_4]}$ to $\M^*(\Sigma_{1,0},\P_3)$. Let $p_1$, $p_2$ and $p_3$ denote the three punctures of the torus, as shown in Figure~\ref{fig:hom}. We denote as $\tau_1$ the right half-twist that interchanges punctures $p_2$ and $p_3$, and as $\tau_2$ the one that interchanges $p_1$ and $p_3$. Note that our numbering of $\tau_i$'s matches the corresponding numbering of flips $\sigma_i$'s according to their action on the punctures. Let also $a_2$, $a_3$ denote the right Dehn twists around the simple closed curves depicted in Figure~\ref{fig:hom}, as before. We define a homomorphism $\Psi\colon{A[F_4]}\to\M^*(\Sigma_{1,0},\P_3)$ on the standard generators as:
\[
\Psi\colon (s_1,s_2,s_3,s_4)\longmapsto (a_3,a_2,\tau_2,\tau_1).
\]

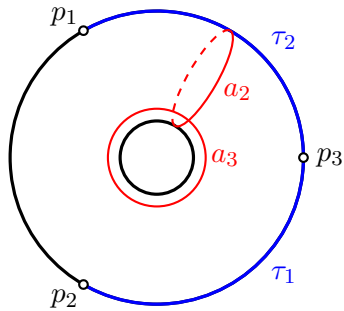
\begin{figure}
\begin{center}
\begin{tikzpicture}[scale=0.65]
\draw [very thick] (0,0) circle (0.75cm);
\draw [black,very thick] (0,0) circle (3cm);
\draw [blue,very thick] (-1.5,2.598) arc [start angle = 120, end angle=0, radius=3cm]; 
\draw [blue,very thick] (3,0) arc [start angle = 0, end angle=-120, radius=3cm]; 
\draw [blue] (2.6,2.4) node {$\tau_2$};
\draw [blue] (2.6,-2.4) node {$\tau_1$};

\color{red}
\draw[black,thick,fill=white,radius=2.5pt] (3,0) circle node [right=1pt] {$p_3$};
\begin{scope}[rotate=-120]
\draw[thick] (-3,0) arc (180:0:1.125cm and 0.30cm);
\draw[thick,dashed] (-3,0) arc (180:360:1.125cm and 0.30cm);
\draw[black,thick,fill=white,radius=2.5pt] (3,0) circle node [below left=-3pt] {$p_2$};
\end{scope}
\begin{scope}[rotate=120]
\draw [black,thick,fill=white,radius=2.5pt] (3,0) circle node [above left=-3pt] {$p_1$};
\end{scope}
\draw (1.65,1.35) node {$a_2$};
\draw [thick] (0,0) circle (1cm);
\draw (1.4,0) node {$a_3$};

\end{tikzpicture}
\caption{The images of the standard generators $(s_1,s_2,s_3,s_4)$ under the homomorphism $\Psi$ from Example~\ref{ex:13}: right Dehn twists $a_3$, $a_2$, and right half-twists $\tau_2$,~$\tau_1$.
\label{fig:hom}}
\end{center}
\end{figure}

\noindent
The properties of Dehn twists and half-twists discussed above guarantee that $\Psi$ is a homomorphism. To prove that $\Psi$ sends the center of $A[F_4]$ to the trivial mapping class, we need to express $\tau_1$ and $\tau_2$ through the standard generators of $\M^*(\Sigma_{1,0},\P_3)$ given in Theorem~\ref{thm:comd4}. We have the following equalities:
\begin{equation}\label{eq:tausigma}\tag{**}
\begin{aligned}
\tau_1 & = \sigma_1\cdot\Delta(a_1,a_3,a_2),\\
\tau_2 & = \sigma_2\cdot\Delta(a_1,a_3,a_4),
\end{aligned}
\end{equation}
where $\Delta(a_1,a_3,a_2)=a_1a_3a_2\cdot a_1a_3\cdot a_1$ is (the image in $\PMod(\Sigma_{1,0},\P_3)$ of) the Garside element of the Artin subgroup of type $A[A_3]\subset A[D_4]$ on the generators $a_1,a_3,a_2$, and, similarly, $\Delta(a_1,a_3,a_4)=a_1a_3a_4\cdot a_1a_3\cdot a_1$. These equalities can be proven by considering the action of $\sigma_i\tau_i$ ($i=1,2$) and the corresponding Garside elements on the curves underlying the standard Dehn twists $a_1,\dots,a_4$ from Figure~\ref{fig:torus}, and observing that these curves cut the surface into punctured disks (\cite[Prop.\,2.8]{FM}).

Now we can treat the right-hand parts of \eqref{eq:tausigma} as words in $A[D_4]\rtimes\S_3$ and compute with the help of the  package \textsf{CHEVIE} that  $\Psi$ sends the generator $(s_1s_2s_3s_4)^6$ of the center of $A[F_4]$ to the element $\delta^7$, where $\delta=(a_1a_2a_3a_4)^3$ is the standard generator of the center of $A[D_4]$. Therefore $\Psi$ defines a homomorphism $\ov{A[F_4]}\longrightarrow\ov{A[D_4]}\rtimes\S_3=\M(\Sigma_{1,0},\P_3)\subset\M^*(\Sigma_{1,0},\P_3)$.

We find with the help of \textsf{GAP} that the projection of $\Psi$ to $\ov{W[D_4]}\rtimes\S_3$ can be conjugated by the element $\hat a_1\hat a_3\hat a_4\hat a_2\hat a_3\cdot \sigma_1\sigma_2\sigma_1$ into the homomorphism $\psi\colon (s_1,s_2,s_3,s_4)\longmapsto (\hat a_3,\hat a_2,\sigma_1,\sigma_2)$, and furthermore, that the image of $\Psi$ has index $9$ in $\M(\Sigma_{1,0},\P_3)$. Thus the reasoning in the last part of the proof of Proposition~\ref{prop:d4f4} applies to $\Psi$ and shows that $\Psi$ is not injective. A straightforward computation shows that the element $(s_1s_2s_3)^6\in\ov{A[F_4]}$ lies in the kernel of $\Psi$.
\end{exa}

Inspired by Example~\ref{ex:13}, we finish this section with the following open question. Recall that we asked in Question~\ref{que:emb} for which Artin groups $A[\Gamma]$ their central quotients $\ov{A[\Gamma]}$ admit embeddings into mapping class groups of punctured surfaces. From Proposition~\ref{prop:d4f4} we deduce that there does not exist an embedding of ${}\ov{A[F_4]}$ into $\M^*(\Sigma_{1,0},\P_3)$ such that its image has finite index in $\M^*(\Sigma_{1,0},\P_3)$. Unlike the case with ${}\ov{A[H_4]}$ (see Corollary~\ref{cor:h4mcg}), this leaves open the question of existence of an embedding of $\ov{A[F_4]}$ into $\M^*(\Sigma_{1,0},\P_3)$ which has the image of infinite index.
\begin{que}
Is there an embedding of $\ov{A[F_4]}$ into $\M^*(\Sigma_{1,0},\P_3)$ such that the image has infinite index in $\M^*(\Sigma_{1,0},\P_3)$?
\end{que}

The analysis in the proof of Proposition~\ref{prop:d4f4} shows that without loss of generality (i.e.\ up to conjugation in $\M^*(\Sigma_{1,0},\P_3)$ and up to the graph automorphisms of $\ov{A[F_4]}$) we can assume that such an embedding, if it exists, has the projection onto $\ov{W[D_4]}\rtimes\S_3$ given on the standard generators by the mapping $(s_1,s_2,s_3,s_4)\longmapsto(\hat a_3,\hat a_2,\sigma_1,\sigma_2)$.

\end{document}